\documentclass[a4paper]{amsart}
\usepackage{amssymb}
\usepackage{amsfonts}
\usepackage{amscd}
\usepackage{verbatim}
\numberwithin{equation}{section} \numberwithin{figure}{section}

\newcommand{\h}{\ensuremath{\mathcal{H}}}

\newcommand{\X}{\ensuremath{\mathcal{X}}}

\newcommand{\Z}{\ensuremath{\mathbb{Z}}}
\newcommand{\p}{\ensuremath{\mathbb{P}}}

\newcommand{\ra}{\ensuremath{\rightarrow}}

\DeclareMathOperator{\m}{mod} \DeclareMathOperator{\Kdim}{Kdim}
\DeclareMathOperator{\mo}{Mod}
\DeclareMathOperator{\gd}{gl.dim}
 
\DeclareMathOperator{\HEIGHT}{ht}
\DeclareMathOperator{\height}{ht} \DeclareMathOperator{\op}{opp}
 
\DeclareMathOperator{\QMod}{QMod}
\DeclareMathOperator{\qmod}{qmod}

\let\cal\mathcal
\def\Ascr{{\cal A}}
\def\Bscr{{\cal B}}
\def\Cscr{{\cal C}}
\def\Dscr{{\cal D}}
\def\Escr{{\cal E}}

\def\Hscr{{\cal H}}

\def\Oscr{{\cal O}}

\def\Sscr{{\cal S}}

\def\Xscr{{\cal X}}

\let\blb\mathbb

\def \PP{{\blb P}}

\def \ZZ{{\blb Z}}

\def \NN{{\blb N}}

\def\id{\text{id}}

\def\Mod{\operatorname{Mod}}
\def\mod{\operatorname{mod}}
\def\Gr{\operatorname{Gr}}

\def\QGr{\operatorname{QGr}}
\def\qgr{\operatorname{qgr}}

\def\gr{\operatorname{gr}}

\def\coh{\mathop{\text{\upshape{coh}}}}

\def\rad{\operatorname {rad}}
\def\gr{\operatorname {gr}}
\def\Spec{\operatorname {Spec}}

\def\Ext{\operatorname {Ext}}
\def\Hom{\operatorname {Hom}}
\def\End{\operatorname {End}}
\def\RHom{\operatorname {RHom}}

\def\uRHom{\operatorname {R\mathcal{H}\mathit{om}}}

\def\coker{\operatorname {coker}}
\def\ker{\operatorname {ker}}

\def\End{\operatorname {End}}

\def\add{\operatorname {add}}
\def\rk{\operatorname {rk}}

\def\r{\rightarrow}

\DeclareMathOperator{\tors}{tors}

\DeclareMathOperator{\HEnd}{\mathcal{E}\mathit{nd}}

\newtheorem{lemma}{Lemma}[section]
\newtheorem{proposition}[lemma]{Proposition}
\newtheorem{theorem}[lemma]{Theorem}
\newtheorem{corollary}[lemma]{Corollary}

\theoremstyle{definition}

{

}

\theoremstyle{remark}

\def\dirlim{\mathop{\vtop{\baselineskip -100pt\lineskip -1pt\lineskiplimit 0pt
\setbox0\hbox{lim}\copy0\hbox to \wd0{\rightarrowfill}}}\limits}
\begin{document}
\title[Grothendieck groups and tilting objects]
{Grothendieck groups and tilting objects}
\author[I. Reiten and M. Van den Bergh]{
Idun Reiten \\ Department of Mathematical Sciences \\ Norwegian
University of Science and Technology \\ 7491 Trondheim \\ Norway
\\ \vspace{0.5cm} Michel Van den Bergh\\ Department WNI, Limburgs
Universitair Centrum \\ Universitaire Campus, Building D
\\3590 Diepenbeek \\ Belgium}

\maketitle

\begin{center}
{\em dedicated to Klaus Roggenkamp for his 60th birthday}
\end{center}

\begin{abstract}
Let $\Cscr$ be a connected noetherian hereditary abelian
 $\Ext$-finite category with Serre functor over an algebraically
closed field $k$, with finite dimensional homomorphism and
extension spaces. Using the classification of such categories from
\cite{RV2}, we prove that if $\Cscr$ has some object of infinite
length, then the Grothendieck group of $\Cscr$ is finitely
generated if and only if $\Cscr$ has a tilting object.
\end{abstract}

\section*{Introduction}
Let $k$ be an algebraically closed field and $\Cscr$ a hereditary
abelian $\Ext$-finite $k$-category. That $\Cscr$ is hereditary
means that $\Ext^{i}(-,-)$ vanishes for $i\ge 2$, and we say that $\Cscr$ is
$\Ext$-finite if $\Ext^{i}(A,B)$ is finite dimensional over $k$
for all $A,B$ in $\Cscr$ and all $i$. A central problem in the
representation theory of artin algebras is to describe such
$\Cscr$ which have a tilting object. This is important in
connection with the investigation of quasitilted algebras, as
introduced in \cite{HRS}.

When a hereditary abelian $\Ext$-finite $k$-category $\Cscr$  has
a tilting object, it is a consequence that the Grothendieck group
$K_0(\mathcal{C})$ is free abelian of finite rank \cite{HRS}. This
suggests the problem of describing the $\Cscr$ for which
$K_0(\Cscr)$ is finitely generated (or free abelian of finite
rank), and to decide to which extent having a finitely generated
Grothendieck group implies the existence of a tilting object.
Relating the existence of a tilting object to properties of the
more widely known notion of Grothendieck group provides a better
insight into  the meaning of the condition of the existence of a
tilting object. In particular, it is interesting to understand in
terms of Grothendieck groups the special role the category
$\coh\cal{X}$ of coherent sheaves on a weighted projective line
\cite{GL1} plays within the larger class of quotient categories of
finitely generated graded modules over commutative noetherian
isolated singularities of Krull dimension two.

We will use the general classification results from \cite{RV2} to
solve the above problems under the additional hypotheses that
$\Cscr$ is noetherian and has a Serre functor (see section
\ref{sec1}). The latter hypothesis is natural since it is a
consequence of the existence of a tilting object. Both additional
hypotheses are satisfied for the quotient categories mentioned
above. Since the properties of $\Cscr$ having a tilting object, a
finitely generated Grothendieck group or a Serre functor are
preserved under derived equivalence of hereditary categories
\cite{HR1,RV2} our results apply more generally to the $\Cscr$
with the additional hypothesis of having a Serre functor and being
derived equivalent to a noetherian hereditary category.

Some of the results on Grothendieck groups of quotient categories
proved in this paper are inspired by similar results for
two-dimensional complete noetherian rings, used as a tool for
classifying maximal orders of finite representation type in
\cite{RV1}.

Section \ref{sec1} is devoted to discussing background material
from various sources, collected together for the benefit of the
reader. In section \ref{sec2} we describe the $\Cscr$ with
finitely generated Grothendieck group. A new criterion for an
object to be a tilting object is given in section \ref{sec3},
which it is interesting also in its own right. In section
\ref{uuuref:5b}  we construct exceptional collections of modules
over a hereditary order over a discrete valuation ring. This is
used in section \ref{sec5}, along with the criterion from section
\ref{sec3}, to construct a tilting object in the category $\coh
\cal{O}$ of coherent  modules over a sheaf $\cal{O}$ of hereditary
orders over $\p^1$. In section \ref{sec6} we give our main result
on the connection between the existence of a tilting object and
the Grothendieck group being finitely generated. Under our
assumptions the conditions turn out to be equivalent if $\Cscr$ is
connected and has some object of infinite length. In section
\ref{sec7} we give some examples  and comments.

We give an appendix proving directly the relationship between the
category of coherent sheaves on a weighted projective line and the
above mentioned category $\coh \cal{O}$.

Hereditary abelian $k$-categories which are $\Ext$-finite,
noetherian and have a Serre functor were classified in \cite{RV2}.
It is also of interest to investigate hereditary abelian
categories which do not satisfy the additional assumptions, for
example with respect to when the Grothendieck group is finitely
generated. In Appendix B we give some sources of examples of
hereditary abelian categories.

\section{Background}
\label{sec1} In this section we provide some background material
from various sources, to provide a better understanding of how our
work fits in.

\subsection{Tilting objects}
Let $\Ascr$ be a triangulated category with the property that
between two objects only a finite number of $\Ext$ are non-zero.
If $T\in\cal{A}$, then add(T) is by definition the smallest
additive category containing $T$ which is closed under finite
direct sums and summands. We say that $T$ is a tilting object if
 $\Ext^{i}_{\cal{A}}(T,T)=0$ for $i \neq 0$ and $\add T$ generates $\cal{A}$
(in the sense that $\cal{A}$ is the smallest subcategory of
$\cal{A}$ containing $\add T$ which is closed under shifts and
cones).

If $\Cscr$ is an abelian category of finite homological dimension
then $T\in\Cscr$ is a tilting object if it is a tilting object in
$D^b(\mathcal{C})$. This definition of a tilting object in $\Cscr$
is equivalent to the usual notion of tilting module (of finite
projective dimension) when $\Cscr$ is the category $\mod \Lambda$
of finitely generated modules for an artin algebra of finite
global dimension, as is seen directly or by using \cite{CPS,Ric}.
For an $\Ext$-finite hereditary abelian $k$-category $\Cscr$ it is
equivalent to the following definition used in \cite{HR2}
(reformulating conditions from \cite{HRS,H2}):
 An
object $T$ in $\Cscr$ is a tilting object if $\Ext^1(T,T)=0$ and
if $\Hom(T,X)=0=\Ext^1(T,X)$ implies $X=0$. In general the
definition in \cite{HRS} is modelled on the definition of a
tilting module of projective dimension at most one, and is hence
different from ours when $\Cscr$ is not hereditary.

When $T$ is a tilting object in $\mod\Lambda$ for an artin algebra
$\Lambda$ of finite global dimension,
 there is an induced  equivalence $D^b(\mod
\Lambda)\to D^b(\mod \End(T)^{\op})$ between bounded derived
categories, and similarly if $T$ is a tilting object in the
category $\coh X$ of coherent sheaves on a smooth projective
variety $X$ \cite{H1,Ric,Baer}. This can easily be extended to the
case of $\coh \cal{O}$ where $\cal{O}$ is a coherent
$\cal{O}_X$-algebra locally of finite global dimension. The
analogous result for tilting objects in $\Ext$-finite hereditary
abelian $k$-categories is given in \cite{HRS}.

In general if $T$ is a tilting object in a triangulated category
$\Ascr$ one may expect an equivalence $\Ascr \cong D^b (\End
(T)^{\op})$. There is recent work in this direction by Keller,
using $A_\infty$-categories, and building on \cite{KB}. In
particular it follows from his results that if $T$ is a tilting
object in an $\Ext$-finite abelian $k$-category $\Cscr$ of finite
homological dimension, then there is an equivalence between
$D^b(\Cscr)$ and $D^b(\End(T)^{\op})$.

It follows from this derived equivalence that when $\Cscr$ is an
$\Ext$-finite abelian $k$-category of finite homological dimension
having a tilting object $T$ , the Grothendieck group  $K_o(\Cscr)$
is isomorphic to $\mathbb{Z}^n$, where $n$ is the number of
nonisomorphic summands of $T$ (see \cite{HRS} for the hereditary
case). Given $T$ in $\mathcal{C}$ with $\Ext^i(T,T)=0$ for $i>0$,
it is an important problem, open even for artin algebras, whether
$T$ having $n$ nonisomorphic summands is sufficient for $T$ to be
a tilting object. It is known in the case of artin algebras when
the projective dimension of $T$ is at most one \cite{Bong}, and
for $\Cscr$ $\Ext$-finite hereditary if $\Cscr$ has some tilting
object \cite{H2}.

\subsection{Hereditary categories and weighted projective lines}
Hereditary abelian $\Ext$-finite categories are of special
interest in connection with quasitilted algebras, as introduced in
\cite{HRS}. The quasitilted algebras are by definition the endomorphism
algebras End$(T)^{\op}$ when $T$ is a tilting object in a
hereditary category. Equivalently, an algebra $\Lambda$ is
quasitilted if and only if gl.dim.$\Lambda\le 2$ and each
indecomposable $X$ in $\mod\Lambda$ has projective or injective
dimension at most one \cite{HRS}. Main examples of categories
$\Cscr$ are $\mod \Lambda$ where $\Lambda$ is a hereditary artin
algebra and $\coh X$ when $X$ is a smooth projective curve. More
generally there are the coherent sheaves on weighted projective
lines \cite{GL1}, which we discuss next.

Let $\lambda = (\lambda_1,\lambda_2,\ldots,\lambda_n)$, $n\geq 2$,
be a finite number of points in $\p^1$, with $\lambda_1=0$,
$\lambda_2=\infty$, $\lambda_3=1$. For a sequence
$e=(e_1,\ldots,e_n)$ of positive integers consider the associated
ring
\begin{equation}
\label{Rdescription}
R=
k[x_1,\ldots,x_n]/[(x_i^{e_i}-x_2^{e_2}+\lambda_ix_1^{e_1})]_{i\geq3}.
\end{equation}
Let $H$ be the abelian group generated by $h_1, \ldots, h_n$ and
with relations $e_1 h_1= \cdots = e_n h_n$. Then $H$ is isomorphic
to $\Z \oplus G$, where $G$ is a finite group \cite{GL1}. $R$ is
then a $H$-graded ring, and  Geigle and Lenzing write $\coh\Xscr$
for the hereditary category gr$_H R/\text{finite length}$, where
$\gr _H R$ denotes the category of finitely generated $H$-graded
$R$-modules with degree zero homomorphisms.  Geometrically $\coh
\Xscr$ can be viewed as the coherent sheaves over a (hypothetical)
space $\Xscr$ which is a generalization of $\mathbb{P}^1$. Hence
Geigle and Lenzing call $\Xscr$ a ``weighted projective line''.

The category $\coh\Xscr$ is a noetherian hereditary abelian category with finite
dimensional homomorphism and extension spaces, which has a tilting
object $T$ such that End$_\h(T)^{\op}$ is a canonical algebra in the
sense of Ringel \cite{Rin}, and actually all canonical algebras are
obtained this way. Like for tilting for finite dimensional algebras
\cite{H1}, there is also induced an equivalence of derived categories
$D^b(\coh \X) \rightarrow \mbox{$D^b($ End$_\h(T)^{\op})$}$ \cite{GL1}. This
setup is used to give an alternative approach to the study of the
module theory for canonical algebras, by first investigating the
hereditary category $\coh \X$. The rings described by
\eqref{Rdescription} are $H$-graded factorial and it is shown in
\cite{Ku} that this property characterizes those rings amongst
the two-dimensional rings.

There are two main known sources of connected hereditary
categories $\Cscr$ with tilting object; the module categories of
finite dimensional hereditary $k$-algebras and the categories
$\coh \X$. In addition there are the hereditary abelian categories
derived equivalent to them. It is conjectured that there are no
more, and in fact this is proved in \cite{L2} for noetherian
hereditary categories and more generally in \cite{HR1,HR2} under
the assumption that $\Cscr$ has at least one nonzero object of
finite length, or at least one directing object, that is, an
object which does not lie on a cycle of nonzero nonisomorphisms.

\subsection{Noetherian hereditary categories with Serre functor}
\label{newssec13}
We recall some essential features of the classification of
noetherian hereditary abelian $\Ext$-finite $k$-categories with
Serre functor. For further details we refer to \cite{RV2}.

Assume that $\Cscr$ is $\Ext$-finite.
A Serre functor for $\Cscr$ is an auto-equivalence $F:D^b(\Cscr)\to
D^b(\Cscr)$ where $D^b(\Cscr)$ denotes the bounded derived
category, such that there are isomorphisms $ \Hom(A,B)
\xrightarrow{\cong}\Hom(B,FA)^\ast$ natural in $A$ and $B$ ($(-)^\ast$
is the $k$-dual). This clearly implies that $\Cscr$ has finite
homological dimension.

For a hereditary abelian $\Ext$-finite $k$-category $\Cscr$ the
existence of a Serre functor implies the existence of almost split
sequences, and the converse holds if $\Cscr$ has no non-zero
projective or injective objects \cite{RV2}.

Let be $\Cscr$ a connected category. It is proved in \cite{RV2}
that if $\Cscr$ is a connected noetherian hereditary $\Ext$-finite
$k$-category with Serre functor then $\Cscr$ has one of the
following forms.
\begin{itemize}
\item[(i)] A category $\coh\Oscr$, where $\Oscr$ is a sheaf of hereditary
orders over a smooth projective curve.
\item[(ii)] $\Cscr$ is
$\mod\Lambda$ for a finite dimensional hereditary $k$-algebra
$\Lambda$.
\item[(iii)] $\Cscr$ is the category of finite dimensional representations
of the quiver $\tilde{A}_n$ with cyclic orientation with
$n<\infty$.
\item[(iv)] $\Cscr$ is derived equivalent to a hereditary category where all
objects have finite length and having an infinite number of nonisomorphic simple objects.
\end{itemize}
The actual result in \cite{RV2} also contains a precise classification
of the categories in (iv).
We also recall from \cite{RV2} that for the categories in (i) there is
an alternative description as follows.

\begin{itemize}
\item[(i')] Categories of the form $\qgr S$ = $\gr S$/finite length,
where $\gr S$ denotes the category of finitely generated graded
modules over a commutative noetherian $\Bbb{Z}$-graded domain
$S=k+S_1+S_2+...+S_i+...$ of Krull dimension two which is finite over
its center, where the $S_i$
are finite dimensional over $k$, and $S$ is an isolated
singularity.
\end{itemize}

The categories $\coh\cal{X}$ have also as we have seen a similar
description, as quotient categories starting with an $H$-graded
ring, where $H$ is not necessarily $\Bbb{Z}$. But it follows from
\cite{GL2}\cite{L1} that one can assume that the rings are
$\Bbb{Z}$-graded, so that the class $\coh\cal{X}$ is a subclass of
the $\qgr S$.

\subsection{Classical hereditary orders.}
\label{ssec13} In this section we collect some well-known
properties of hereditary orders. We will loosely refer to a
classical hereditary order as an order $\Lambda$ in a central
simple algebra $A$ over a field $K$ which is hereditary. Let $R$
be the center of $\Lambda$. According to \cite{RS} $R$ is a
Dedekind ring. (This fact does not seem to be contained in exactly
this form in \cite{Reiner}).

Assume that $R$ is a discrete valuation ring with maximal ideal
$m$. Then according to \cite{Reiner}  the radical $I$ of $\Lambda$
is invertible. Furthermore by \cite{Reiner} there is an integer
$e$ such that $I^e=m\Lambda$, called the ramification index of
$\Lambda/R$. It follows from the structure theory of hereditary
orders in \cite{Reiner} that if $e=1$ then $\Lambda$ is maximal.
If the converse is true then we say that $A/K$ is unramified. This
happens for example if $A=M_n(K)$.

If $R$ is not  a discrete valuation ring then by localizing one
defines ramification indices $e_P$ for the non-zero primes in $R$.
By analyzing  $\Dscr=\operatorname{Hom}(\Lambda,R)$ it follows
easily that $\Lambda/R$ ramifies in only a finite number of
primes.

\section{Finitely generated Grothendieck groups}
\label{sec2} Let $\Cscr$ be a noetherian $\Ext$-finite hereditary
abelian $k$-category with Serre functor, where $k $ is an
arbitrary field. In this section we describe which $\Cscr$ have
finitely generated Grothendieck group. For this we use the
classification theorem from \cite{RV2} in the form recalled in
\S\ref{newssec13}. The main problem we need to deal with is when the category
$\coh\cal{O}$ of coherent modules over a sheaf $\cal{O}$ of
hereditary orders over a smooth projective curve $X$ has finitely
generated Grothendieck group.

Let $X$ be a regular connected curve over a field $k$. Let $K$ be
the function field of $X$ and let $A$ be a central simple algebra
over $K$. Let $\Oscr$ be a sheaf of hereditary orders in $A$ over
$\Oscr_X$. Thus locally $\Oscr$ is a hereditary order over a
Dedekind ring (in the sense of \cite{Reiner}).

If $\Ascr$ is a sheaf of rings on a topological space $Z$ then we
use
 the
notation $\operatorname{coh}(\Ascr)$ for the category of coherent
$\Ascr$-modules and we write $K_0(\Ascr)$ for
$K_0(\operatorname{coh}(\Ascr))$.

Our first aim is to give some results on $K_0(\Oscr)$.
\begin{proposition}
\label{ref:3.1a} Let $\Oscr$ be as above. Let $x_1,\ldots,x_t$ be
the points in which $\Oscr$ ramifies, and
 let
$e_1,\ldots,e_t$ be the corresponding ramification indices (see
 \S\ref{ssec13}).

Put
\begin{equation}
\label{ref:2a} r=\sum_i(e_i-1).
\end{equation}
Let $\bar{\Oscr}$ be a maximal order lying over $\Oscr$.
 Then
\[
K_0(\Oscr)\cong K_0(\bar{\Oscr})\oplus \Z^r.
\]
If $k$ is algebraically closed then
\[
K_0(\Oscr)\cong K_0(\Oscr_X)\oplus
 \Z^r.
\]
\end{proposition}
\begin{proof}
The hereditary orders in $A$ containing $\Oscr$ form a partially
ordered set which we will denote by $\Hscr(\Oscr)$.

If $\Oscr'\in \Hscr(\Oscr)$ lies minimally over $\Oscr$ then one
proves exactly as in \cite[Thm 1.14]{RV1} that
$K_0(\Oscr)=K_0(\Oscr')\oplus \Z$.

Let $\bar{\Oscr}$ be a maximal order lying over $\Oscr$. We deduce
that $K_0(\Oscr)\cong \Z^{r}\oplus K_0(\bar{\Oscr})$, where $r$ is
the length of a maximal chain in $\Hscr(\Oscr)$, starting in
$\Oscr$ and ending in $\bar{\Oscr}$.

A local computation shows that $r$ is given by the formula
\eqref{ref:2a}, which finishes the proof of the computation of
$K_0(\Oscr)$.

If $k$ is algebraically closed then by Tsen's theorem \cite[p.\
374]{C} one has that $\bar{\Oscr}\cong \End_{\Oscr_X}(\Escr)$
where $\Escr$ is a vector bundle of rank $n$ on $X$. Hence by
Morita theory $K_0(\bar{\Oscr})\cong K_0(\Oscr_X)$.
\end{proof}
\begin{corollary}
\label{corcor} Assume $k$ is a algebraically closed. Then
$K_0(\Oscr)$ is
 finitely
generated if and only if $X$ is an open subset of
 $\mathbb{P}^1$.
\end{corollary}
\begin{proof}
  By the previous proposition it suffices to prove this for
  $\Oscr=\Oscr_X$.  Let $\bar{X}$ be the regular projective curve
  associated to the function field  of $X$ \cite{F}. Then $\bar{X}$ is a
  regular compactification of $X$.  In particular $\bar{X}-X$ is a finite number of
  points, whence by the localization sequence $K_0(\Oscr_X)$ is
  finitely generated if and only if $K_0(\Oscr_{\bar{X}})$ is finitely
  generated. Hence we may assume that $X$ is projective. By \cite[Ex.
  II.6.12, Rem. IV.4.10.4]{Har2} one has $K_0(\Oscr_X)\cong \Z^2\oplus J(k)$
  where $J(k)$ denotes the $k$-points of the Jacobian of $X$. It is
  well-known that $J(k)$ is not finitely generated if $X \ncong \mathbb{P}^1$
  (for example because in that case $J(k)$ is non-trivial and
  divisible \cite{M}).
\end{proof}

Combining with \S\ref{newssec13} we now get the following main result of this
section.

\begin{theorem}
\label{newthm23}
Let $\Cscr$ be a connected noetherian $\Ext$-finite hereditary
abelian $k$-category with Serre functor where $ k $ is an
algebraically closed field. Then $K_o(\Cscr)$ is finitely
generated (free abelian) if and only if $\Cscr$ has one of the
following forms.
\begin{enumerate}
\item
$\mod\Lambda$ where $\Lambda$ is an indecomposable finite
dimensional hereditary $k$-algebra.
\item
Finite dimensional representations over $\tilde{A}_n$ with $n$
finite and cyclic orientation.
\item
$\coh\cal{O}$ where $\cal{O}$ is a sheaf of hereditary
$\cal{O}_X$-orders with $X=\Bbb{P}^1$.
\end{enumerate}
\end{theorem}

\begin{proof}
It is clear that the categories in 1. and 2. have finitely
generated (free abelian) Grothendieck groups, and that
$K_o(\Cscr)$ is not finitely generated if $\Cscr$ is derived
equivalent to a hereditary category $\Cscr'$ with all objects of
finite length and an infinite number of nonisomorphic simple
objects. In view of \S\ref{newssec13} the proof is completed by using Corollary
\ref{corcor}.
\end{proof}

\section{A criterion for deciding if an object is a tilting object.}
\label{sec3} The aim of this section is to give a criterion for an
object of projective dimension at most one in an abelian category
to be a tilting object. For this we need to recall some results on
semiorthogonal pairs in triangulated categories from
\cite{Bond,BK}.

Let $\Ascr$ be a triangulated category and let $\Bscr$, $\Cscr$ be
 two
strict ($=$ closed under isomorphisms) full triangulated
subcategories of $\Ascr$.  $(\Bscr,\Cscr)$ is said to be a
\emph{semi-orthogonal
  pair} if $\Hom_\Ascr(B,C)=0$ for $B\in\Bscr$ and $C\in
 \Cscr$.
Define
\[
\Bscr^\perp=\{A\in\Ascr\mid \forall B\in\Bscr
 :\Hom_\Ascr(B,A)=0\}
\]
${}^\perp\Cscr$ is defined similarly.

If $\Sscr$ is a class of objects in $\Ascr$ then the
(triangulated) category generated by $\Sscr$ is the smallest
subcategory of $\Ascr$ which is closed under shifts, cones, and
isomorphisms.

The following result is a slight variation of the statement of
\cite[Lemma 3.3.1]{Bond} (see also  \cite[\S 1]{BK}).
\begin{lemma}
\label{zzzref:4.1a} The following conditions are equivalent for a
  semi-orthogonal pair $(\Bscr,\Cscr)$.
\begin{enumerate}
\item
$\Bscr$ and $\Cscr$ generate $\Ascr$.
\item For every $A\in\Ascr$ there exists a distinguished triangle $B\r
  A\r C$ with $B\in\Bscr$ and $C\in\Cscr$.
\item $\Cscr=\Bscr^\perp$ and  the inclusion functor
  $i_\ast:\Bscr\r\Ascr$ has a right adjoint $i^!:\Ascr\r \Bscr$.
\item $\Bscr={}^\perp\Cscr$ and  the inclusion functor $j_\ast:\Cscr\r
  \Ascr$ has a left adjoint  $j^\ast:\Ascr\r \Cscr$.
\end{enumerate}
If one of these conditions holds then the functors $i^!,j^\ast$
are exact and the triangles in 2. are (for a fixed $A$) unique up
to unique isomorphism. They are necessarily of the form
\begin{equation}
\label{zzzref:5a} i_\ast i^! A\r A \r j_\ast j^\ast A
\end{equation}
where the maps are obtained by adjointness from the identity maps
$i^!A\r i^! A$ and $j^\ast A \r j^\ast A$. In particular triangles
as in 2. are functorial.
\end{lemma}
If any of the conditions of the previous lemma holds then we say
 that
$(\Bscr,\Cscr)$ is a semi-orthogonal decomposition of $\Ascr$.

Now for the rest of this section let $k$ be a field.
 All
categories (abelian or triangulated) will be $k$-linear and have
finite dimensional $\Hom$'s and $\Ext$'s. We assume furthermore
that for any pair $A,B$, there are only a finite number of
 non-zero
$\Ext^i(A,B)$.

Let $\Ascr$ be a triangulated or abelian category. For an object
$T$ in $\Ascr$ we denote by $T^\perp$ the full subcategory of
$\Ascr$ whose objects are the $C$ in $\Ascr$ with
$\Ext^i(T,C)=0$ for all $i$. We say that an
 object
$T\in\Ascr$ is exceptional if $\Ext^i_\Ascr(T,T)=0$ for $i>0$ and
$\End_\Ascr(T)$ is a (finite dimensional) division
 algebra. A sequence of exceptional objects $T_1,\ldots,T_n$ is an
 exceptional collection if  $\Ext^*_\Ascr(T_i,T_j)=0$ for $j>i$. An
 exceptional collection is strongly exceptional if
 $\Ext^t_\Ascr(T_i,T_j)=0$ for $t>0$ and all $i,j$.

The following is proved in  \cite{BK}.
\begin{lemma}
\label{zzzref:4.2a} Assume that $T_1,\ldots,T_n$ is an exceptional
collection in a triangulated category $\Ascr$. Let $\Bscr$ be the
triangulated subcategory of $\Ascr$ generated by $T_1,\ldots,T_n$
and put $T=\oplus_i T_i$. Then $\Ascr$ has a semi-orthogonal
decomposition given by $(\Bscr,\Bscr^{\perp})=(\Bscr,T^{\perp})$.
\end{lemma}
\begin{proof} For the convenience of the reader we repeat the proof.
We have to show that $\Ascr$ is generated by $\Bscr$ and
$T^\perp$.
 Let $\Bscr_1$ be the full subcategory of $\Ascr$ consisting
  of objects isomorphic to finite direct sums of the form $\oplus_j
  T_1[j]^{a_j}$. Then $\Bscr_1$ is a strict triangulated subcategory
  of $\Ascr$. This can be deduced from the fact that the formation of
  triangles in $\Ascr$ is compatible with direct sums \cite[Cor.\
  II.1.2.5]{V}.  Sending $A\in\Ascr$ to $\oplus_i
  \Ext^i(T_1,A)\otimes_D T_1[-i]$ defines a right adjoint to the
  inclusion $\Bscr_1\hookrightarrow \Ascr$. This yields a
  semi-orthogonal decomposition of $\Ascr$ given by
  $(\Bscr_1,\Bscr_1^\perp)=(\Bscr_1,T_1^\perp)$. In particular $\Ascr$
 is generated by $\Bscr_1$ and $T_1^\perp$. Now we repeat this
  construction with $T_2\in T_1^\perp$. So if $\Bscr_2$ is the full
  subcategory of $\Ascr$ consisting of objects isomorphic to finite
  direct sums of the form $\oplus_j T_2[j]^{a_j}$ then we have that
 $T_1^\perp$ is generated by $\Bscr_2$ and $(T_1\oplus
 T_2)^\perp$. Continuing this procedure we find that $\Ascr$ is
 generated by $\Bscr_1,\ldots \Bscr_n$ and $T^\perp$. This finishes
 the proof.
\end{proof}

We point out the following consequence of Lemma \ref{zzzref:4.2a}.
\begin{corollary}
Let $T_1,\cdots,T_n$ be a strongly exceptional collection in a
triangulated category $\Ascr$ satisfying the above assumptions.
Then $T$ is a tilting object if and only if $T^\perp=0$.
\end{corollary}

We shall need that semi-orthogonal decompositions behave nicely
with respect to Grothendieck groups.

\begin{lemma}
\label{zzzref:4.3a} Assume that $(\Bscr,\Cscr)$ is a
semi-orthogonal
  decomposition for a triangulated category $\Ascr$. Then
  $K_0(\Ascr)\cong K_0(\Bscr)\oplus K_0(\Cscr)$.
\end{lemma}
\begin{proof}
The inclusions $\Bscr,\Cscr\subset \Ascr$ define a map
$K_0(\Bscr)\oplus K_0(\Cscr)\r K_0(\Ascr)$. An inverse to this map
is given by sending $[A]$ to $[i^!A]\oplus [j^\ast A]$ (see
 Lemma
\ref{zzzref:4.1a} for notations).
\end{proof}
\begin{lemma} \label{zzzref:4.4a}
Assume that $\Ascr$ is a triangulated category, and let
$T_1,\ldots,T_n\in\Ascr$ be an exceptional collection.
 Put $T=\oplus_i T_i$. Then
  $K_0(\Ascr)\cong \Z^n\oplus
 K_0(T^{\perp})$.
\end{lemma}
\begin{proof}
Using the same method as in Lemma \ref{zzzref:4.2a} we find
inductively using Lemma \ref{zzzref:4.3a} that
$K_0(\Ascr)=\oplus_i K_0(\Bscr_i)\oplus K_0(T^\perp)$. Now it is
easy to see that sending $\oplus_j T_i[j]^{a_j}$ to $\sum_j (-1)^j
a_j$ defines an isomorphism $K_0(\Bscr_i)\cong \ZZ$. This proves
what we want.
\end{proof}

If $\Bscr$ is an abelian category and $B\in\Bscr$ then we will say
that $B$ has projective dimension $\le r$ if $\Ext^s_\Bscr(B,-)=0$
 for
$s>r$.

The above results have a counterpart for abelian categories
provided we work with exceptional objects of projective dimension
$\le 1$. This follows from the following lemma.
\begin{lemma}
\label{zzzref:4.5a} Assume that $\Bscr$ is an abelian category and
let
  $T$ be an  object in $\Bscr$ of projective dimension $\le
  1$.
If an object in $D^b(\Bscr)$ is (right) perpendicular to $T$
  then so is its homology. In
 particular
$D^b_{T^\perp_\Bscr}(\Bscr)=T^\perp_{D^b(\Bscr)}$.
\end{lemma}
\begin{proof}
Let $B\in T^\perp_{D^b(\Bscr)}$ and assume that $n$ is maximal
such that $H^n(B)\neq 0$. Then there is a
 triangle
\begin{equation}
\label{zzzref:6a} \tau_{<n} B\r B\r H^n(B)[-n]\r
\end{equation}
Applying $\Hom(T,-)$ yields
 injections
$\Ext^i(T,H^n(B))\hookrightarrow\Hom(T,(\tau_{<n} B)[n+1+i])$.
 Since for $i\ge 0$ the non-trivial homology of $(\tau_{<n} B)[n+1+i]$
occurs in degrees $\le -2$ and since the projective dimension of
$T$ is less than or equal to 1 it follows that $\Hom(T,(\tau_{<n}
B)[n+1+i])=0$ for $i\ge 0$. In particular $\Ext^i(T,H^n(B))=0$ for
$i\ge 0$. Since trivially
 $
\Ext^i(T,H^n(B))=0$ for $i<0$ it follows that $H^n(B)\in T^\perp$.
 But
then it follows from \eqref{zzzref:6a} that $\tau_{<n} B\in
T^\perp$. Repeating this procedure with  $B$ replaced by
$\tau_{<n}
 B$
eventually yields that the homology of $B$ is in
 $T^\perp$.
\end{proof}
As a corollary one obtains a proof of the following standard
result.
\begin{corollary} Let $T$ and $\Bscr$ be as in the previous
 lemma. Then $T^\perp_\Bscr$ is an abelian category.
\end{corollary}
\begin{proof} If $f:A\r B$ is a map in $T^\perp_\Bscr$ then one has to
  show that $\ker f$, $\coker f\in T^\perp_\Bscr$. Since the complex
  represented by $f$ clearly lies in $T^\perp_{D^b(\Bscr)}$, this
  follows from the previous lemma.
\end{proof}

As a consequence of Lemma \ref{zzzref:4.5a} we obtain the following result on
Grothendieck groups.

\begin{corollary}
\label{zzzref:4.7a} Assume that $\Bscr$ is an abelian category.
 Assume that $T_1,\ldots,T_n\in\Bscr$ is an exceptional collection
 consisting of objects of projective dimension $\le 1$.
 Let $T=\oplus_i T_i$. Then
  $K_0(\Bscr) \cong \ZZ^n\oplus
 K_0(T^{\perp})$.
\end{corollary}
\begin{proof}
By  Lemmas \ref{zzzref:4.4a} and \ref{zzzref:4.5a} one has
\begin{align*}
K_0(\Bscr)&\cong K_0(D^b(\Bscr))\\ &\cong \ZZ^n\oplus
K_0(T^\perp_{D^b(\Bscr)})\\ &\cong \Z^n\oplus
K_0(D^b_{T^\perp_\Bscr}(\Bscr))\\ &\cong \Z^n\oplus
K_0(T^\perp_\Bscr)
\end{align*}
\def\qed{}\end{proof}

We now get the main result of this section.
\begin{corollary}
\label{zzzref:4.8a}
  Assume that $\Bscr$ is an abelian category of finite Krull
  dimension. Let $T=\oplus_{i=1}^n T_i$ be as in Corollary
  \ref{zzzref:4.7a},
  but assume is addition that $\Ext^1(T_i,T_j)=0$ for all $i,j$.
  If $n=\rk K_0(\Bscr)$ then $T$ is a tilting object in
  $\Bscr$.
\end{corollary}
\begin{proof} By Lemma \ref{zzzref:4.2a} we have to show
  $T^{\perp}=0$. By Lemma \ref{zzzref:4.5a} it follows that $T^\perp$ is equal to
  $D^b_{T^\perp_\Bscr}(\Bscr)$. So it is sufficient to show that
  $\Hscr\overset{\text{def}}{=}T^\perp_\Bscr=0$.

  By Corollary \ref{zzzref:4.7a} one has $\rk K_0(\Hscr)=0$. Since
  $\Hscr$ is an abelian subcategory of $\Bscr$, it also has finite Krull
  dimension. In particular if $\Hscr\neq 0$ there is a quotient
  category $\Cscr$ of $\Hscr$ which has finite length. Selecting a
  simple object in $\Cscr$ yields a rank function on $\Hscr$ which
  is non-trivial.  Hence $\rk
  K_0(\Hscr)>0$. This yields a contradiction.
\end{proof}

It would be interesting to know if for a nonzero hereditary
abelian $k$-category $\Bscr$ with finite dimensional homomorphism
and extension spaces we must have $K_0(\Bscr)\neq 0$.

\section{Strongly exceptional collections for hereditary orders over
discrete valuation rings}
\label{uuuref:5b} In order to construct a tilting object in the
category $\coh \cal{O}$ of coherent modules for a hereditary order
$\cal{O}$ over $\p^1$ we need to produce some exceptional collections of
modules for hereditary orders over discrete valuation rings. We start
by recalling some properties for such orders.  For simplicity we
restrict ourselves to hereditary orders contained in a matrix ring
since that is the only case we will need.

Let $R$ be a discrete valuation ring and let $m$ be its maximal
ideal. Furthermore let $K$ be the quotient field of $R$ and put
$A=M_n(K)$ . Let $\Delta$ be a hereditary order in $A$ in the
sense of \cite{Reiner}.

Thus there exist strictly positive integers $n_1,\ldots, n_t$ such
 that
$n=n_1+\cdots +n_t$ and such that $\Delta$ is isomorphic
 to
\begin{equation}
\label{uuuref:7a}
\begin{pmatrix}
R_{n_1\times n_1}& m_{n_1\times n_2}&\cdots &\\ R_{n_2\times n_1}&
R_{n_2\times n_2} &\cdots &\\ \vdots &\vdots &\ddots&\\ &&&
R_{n_t\times n_t}
\end{pmatrix}
\end{equation}
Here $(-)_{a\times b}$ is a shorthand for $M_{a\times b}(-)$.
Strictly speaking this is proved in \cite{Reiner} only in the case
that $R$ is complete, but as is remarked in \cite[bottom
of p.\ 364]{Reiner} the result remains valid in the case we consider.

For $i=1,\ldots,t$ put $p_i=\sum_{j\le i} n_j$, $q_i=\sum_{j>i}
n_j$. Also let $P_i$ be the $i$'th indecomposable projective for
$\Delta$.

Thus by
 definition
\[
P_i=\begin{pmatrix} m_{p_{i-1}}\\ R_{q_{i-1}}
\end{pmatrix}
\]
where $(-)_a$ now stands for $M_{a\times 1}(-)$.

Clearly $P_1\supset P_2\supset \cdots \supset P_t$. Set
$S_i=P_1/P_{i}$ for $i=2,\ldots,t$ and for the same indexes put
 define
\[
\Delta_i=
\begin{pmatrix}
R_{p_{i-1}\times p_{i-1}}& m_{p_{i-1}\times q_{i-1}}\\
R_{q_{i-1}\times p_{i-1}}& R_{q_{i-1}\times q_{i-1}}
\end{pmatrix}
\]
The $\Delta_i$ are  submaximal
orders containing $\Delta$ and  $S_i$ is a simple
$\Delta_i$-module.

Using this observation one proves:
\begin{lemma}
One has
\[
\Delta_i\otimes_{\Delta}S_j=
\begin{cases}
0&\text{if $i>j$}\\ S_i&\text{if $i\le j$}
\end{cases}
\]
\end{lemma}
\begin{proof}
Since one has $\Delta_i\otimes_\Delta P_1=\Delta_i\otimes_\Delta
 \Delta_i\otimes_{\Delta_i}P_1=\Delta_i\otimes_{\Delta_i} P_1=P_1$
 it
follows that
\[
\Delta_i\otimes_{\Delta} S_j=P_1/\Delta_i P_j
\]
It is now easy to see that
\[
\Delta_iP_j=
\begin{cases}
P_1&\text{if $i>j$}\\ P_i&\text{if $i\le j$}
\end{cases}
\]
which yields the
 result.
\end{proof}
\begin{corollary}
\label{uuuref:5.2a} $\Ext^1_\Delta(S_j,S_i)=0$,
 and
\[
\Hom_\Delta(S_j,S_i)=
\begin{cases}
0&\text{if $i>j$}\\ R/m&\text{if $i\le j$}
\end{cases}
\]

\end{corollary}
\begin{proof} One
 has
\[
\Ext^p_\Delta(S_j,S_i)=\Ext^p_{\Delta_i}(\Delta_i\otimes_\Delta
S_j,S_i)
=
\begin{cases}
0&\text{if $i>j$}\\ \Ext^p_{\Delta_i}(S_i,S_i)&\text{if $i\le j$}
\end{cases}
\]
It is easy to see that $\Hom(S_i,S_i)=R/m$ and one  verifies
directly that
$
\Ext^1_{\Delta_i}(S_i,S_i)=0 $.
\end{proof}

\begin{corollary}
$(S_i)_{i=2,\ldots, t}$ is a strongly exceptional collection for
$\Delta$.
\end{corollary}

One also has
\begin{lemma}
\label{uuuref:5.3a} $\Ext^p_\Delta(S_j,P_1)=0$ for all $p$.
\end{lemma}
\begin{proof}
Put $\Gamma=M_n(R)$. Then one has
$
\Ext^p_\Delta(S_j,P_1)=\Ext_{\Gamma}^p(\Gamma\otimes_\Delta
 S_j,P_1)=0
$
\end{proof}

\section{Existence of tilting objects for hereditary orders on
 $\PP^1$}
\label{sec5} In this section $X$ will be $\PP^1$ for an
algebraically field $k$. Let $K$ be the function field of $X$ and
let $A=M_n(K)$.
 Let $\Oscr$ be a sheaf of hereditary orders
in $A$. Thus locally $\Oscr$ is a hereditary order over a Dedekind
ring  (in the sense of \cite{Reiner}).

To compute global $\Ext$'s in $\coh(\Oscr)$ below we will use the
fact that \cite[Prop.
 II.5.3]{G}
\begin{equation}
\label{ref:8a} \RHom_\Oscr(A,B)=R\Gamma(X,\uRHom_\Oscr(A,B))
\end{equation}
together with the fact that $\uRHom_\Oscr(A,B)$ can be computed
locally. That is, if $x\in X$ then
$\uRHom_\Oscr(A,B)_x=\RHom(A_x,B_x)$.

Let $\bar{\Oscr}$ be a maximal order in $A$ lying over $\Oscr$. By
Tsen's theorem there exists a vector bundle $\Escr$ of rank $n$ on
$X$ such that $\bar{\Oscr}=\HEnd_{\Oscr_X}(\Escr)$.

 Fix $i$ and
put $R_i=\Oscr_{X,x_i}$, $\Delta_i=\Oscr_{x_i}$. Then $\Delta_i$
is isomorphic to  an order of the form \eqref{uuuref:7a} with
$t=e_i$. We choose this isomorphism in such a way that it extends
to an isomorphism between $\bar{\Oscr}_{x_i}$ and $M_n(R_i)$. Let
us write $S_{ij}$ ($i=1,\ldots t,j=2,\ldots,e_i$) for
  the finite length $\Delta_i$-modules which were denoted by $S_j$ in
 Section \ref{uuuref:5b}.
 We consider the $S_{ij}$ as $\Oscr$-modules.

Since $\coh \Oscr_X$ has a tilting object, for example given by
$\Oscr_X\oplus \Oscr_X(-1)$, the same is true for $\bar{\Oscr}$ by
Morita theory. We will take $\bar{T}=\Escr\oplus \Escr(-1)$ as a
tilting object in $\coh(\bar{\Oscr})$. By the choice of the local
isomorphisms, $\Escr_{x_i}$ will correspond to the projective
$\Delta$-module denoted by $P_1$ in the previous
 section.

\begin{proposition}
\label{newprop51}
$T=  \bigoplus_{ij}S_{ij}\oplus \Escr\oplus \Escr(-1) $ is a
tilting object in $\coh(\Oscr)$.
\end{proposition}
\begin{proof} Since $K_0(\Oscr_{\PP^1})\cong \Z^2$ it follows from Proposition
  \ref{ref:3.1a} that the number of summands of $T$ is equal to the rank
  of $K_0(\Oscr)$.

  We want to show that the summands of $T$ are a strongly exceptional collection. To do this we
  have to compute the $\Ext^*(-,-)$ between the summands of $T$. We
  first compute the $\uRHom$'s using Corollary \ref{uuuref:5.2a} and
  lemma \ref{uuuref:5.3a}.  The result is as follows.
\begin{equation}
\begin{array}{|c|c|c|c|}
\hline & S_{kl} & \Escr & \Escr(-1)\\ \hline S_{ij}   & \ast & 0 &
0\\ \hline \Escr & \Oscr_{x_k} &  \Oscr_X & \Oscr_X(-1)\\ \hline
\Escr(-1) &  \Oscr_{x_k}& \Oscr_X(1) & \Oscr_X\\ \hline
\end{array}
\end{equation}
For the square marked `$\ast$' we have (using Corollary
\ref{uuuref:5.2a})
\[
\uRHom(S_{ij},S_{kl})=
\begin{cases} 0&\text{if $i\neq k$}\\
0&\text{if $i=k$ and $l>j$}\\ \Oscr_{x_i}&\text{otherwise}
\end{cases}
\]
It now follows immediately from \eqref{ref:8a} that $\Ext^\ast$ is
zero between the summands of $T$. For the $\Hom$'s we find:
\begin{equation}
\begin{array}{|c|c|c|c|}
\hline & S_{kl} & \Escr & \Escr(-1)\\ \hline S_{ij}   & \ast & 0 &
0\\ \hline \Escr & k &  k & 0\\ \hline \Escr(-1) &  k& k^2 & k\\
\hline
\end{array}
\end{equation}
with the `$\ast$' entry given by
\[
\Hom(S_{ij},S_{kl})=
\begin{cases} 0&\text{if $i\neq k$}\\
0&\text{if $i=k$ and $l>j$}\\ k&\text{otherwise}
\end{cases}
\]
So it follows in particular that $T$ is defined by a strongly exceptional
collection. We are now done by Corolllary \ref{zzzref:4.8a}.
\end{proof}
\section{Finitely generated Grothendieck groups and existence of tilting objects.}
\label{sec6} In this section we combine our previous results to
get our desired connection between existence of tilting objects
and the Grothendieck group being finitely generated.

The main result of this paper is the following.
\begin{theorem}
Let $\Cscr$ be a connected noetherian hereditary abelian
 $\Ext$-finite $k$-category with Serre functor, where $k$ is an
algebraically closed field. Then the following are equivalent.
\begin{itemize}
\item[(a)] $K_o(\Cscr)$ is finitely generated.
\item[(b)] \begin{itemize} \item[(i)] $\Cscr$ has a tilting object or
\item[(ii)] $\Cscr$ is the category of finite dimensional
representations of the quiver
 $\tilde{A}_n$ with cyclic orientation for some $n<\infty$.
\end{itemize}\end{itemize}
\end{theorem}

\begin{proof}
$(b)\Rightarrow(a)$. We have already pointed out that $(b)(i)$
implies $(a)$ \cite[I.4.6]{HRS}, and $(b)(ii)$ implies $(a)$ is obvious.

$(a)\Rightarrow(b)$. Assume that $K_o(\Cscr)$ is finitely
generated. If $\Cscr=\mod\Lambda$ for a finite dimensional
hereditary $k$-algebra $\Lambda$, then $\Cscr$ has a tilting
object. If $\Cscr=\coh\Oscr$ where $\Oscr$ is a sheaf of
hereditary orders over $\PP^1 $, it follows from Proposition \ref{newprop51}
that $\Cscr$ has a tilting object. Hence we are done using Theorem
\ref{newthm23}.
\end{proof}

Actually, the following related result is also of interest.
\begin{theorem}
\label{newthm62}
\label{yetanothermainresult} Let $\Cscr$ be a connected noetherian
$\Ext$-finite
  hereditary category which has no projectives or injectives and which
  has an object which is not of finite length. Then the following are
  equivalent.
\begin{enumerate}
\item \label{itemc} $\Cscr$ has a tilting object.
\item  \label{itemb}
$\Cscr$ is derived equivalent to a finite dimensional algebra.
\item  \label{itema}
$\Cscr$ has almost split sequences and $K_0(\Cscr)$ is finitely
  generated.
\item  \label{itemd}
$\Cscr$ is of the form $\coh(\Oscr)$ where $\Oscr$ is a sheaf of
  hereditary $\Oscr_{\PP^1}$-orders.
\item  \label{iteme}
$\Cscr$ is of the form $\coh\Xscr$ for a weighted projective
  line $\Xscr$.
\end{enumerate}
\end{theorem}

\begin{proof}
$1.\Rightarrow 2$. When $\Cscr$ has a tilting object $T$, it
follows from \cite[I, Th. 4.6]{HRS} that $\Cscr$ is derived
equivalent to the finite dimensional algebra End$(T)^{\op}$.

$2.\Rightarrow 3$. Since the hereditary category $\Cscr$ is
derived equivalent to a finite dimensional algebra $\Lambda$, it
follows that $\Lambda$ must have finite global dimension. Hence
$\mod\Lambda$, and consequently $\Cscr$, has a Serre functor
\cite{H1}. Then it follows that $\Cscr$ has almost split sequences
\cite{RV2}.

Since $K_0(\mod\Lambda)$ is finitely generated, it follows that
$K_o(\Cscr)$ is finitely generated because this property is an
invariant of derived equivalence.

$3.\Rightarrow 4$. Since $\Cscr$ has almost split sequences and no
nonzero projectives or injectives, it follows that $\Cscr$ has a
Serre functor \cite{RV2}. Since $\Cscr$ has some object of
infinite length, it follows from Theorem \ref{newthm23} that $\Cscr$ is of
the form $\coh(\Oscr)$ where $\Oscr$ is a sheaf of hereditary
$\Oscr_{\p^{1}}$-orders.

$4.\Rightarrow 1$. This follows from Proposition \ref{newprop51}.

$5.\Leftrightarrow 1$. That $\coh\Xscr$ for a weighted projective
line $\Xscr$ has a tilting object follows from \cite{GL1}, and
$1.\Rightarrow 5$. follows from \cite{L2}.
\end{proof}
In an appendix we give
for completeness an independent proof of $4.\Leftrightarrow 5$,
hence providing a proof of Theorem \ref{newthm62} without using \cite{L2}.

\section{Examples and comments}

\label{sec7}

In this section we give some examples and comments without proofs,
related to the material in this paper.

We start by pointing out how to obtain some concrete examples of
categories $\qgr S$. Translation quivers $\ZZ \Delta$, where
$\Delta$ is an extended Dynkin diagram occur as AR-quivers for the
graded reflexive modules over invariant rings $S=k[X,Y]^G$ where
$G$ is a finite group and $k$ is an algebraically closed field of
characteristic zero (see \cite{ARS} for the definition of
AR-quiver). The corresponding mesh category for $\ZZ \Delta$ is
then a full subcategory of $\qgr S$ whose objects have no nonzero
summands of finite length. We obtain a (finite) basis for
$K_0(\qgr S)$ by considering vertices given by a "section".
Actually such a set of vertices corresponds to a tilting object.
For two-dimensional $Z$-graded rings $S'$ of finite (graded)
representation type we have that $\qgr S'$ is equivalent to some
$\qgr k[X,Y]^G$. The rings $k[X,Y]$ are Gorenstein. It is however
not true for a two-dimensional isolated singularity in general
that there is a commutative Gorenstein ring $S$ with $\qgr S'$
equivalent to $\qgr S$(see \cite{L1}).

While there is a lot of analogy with the work in \cite{RV1}, we
note that there are also some differences. In the complete case,
for the rings $\Lambda$ of finite representation type, considered
as orders, the rank $n$ of $K_0(q(\Lambda))$ gives information on
how far $\Lambda$ is from being a maximal order (of finite
representation type). In this case there was a chain
$\Lambda=\Lambda_1 \subset \cdots \subset \Lambda_n$ of orders
with $\Lambda_n$ maximal and such that there is no refinement of
the chain. Then $K_0(q(\Lambda))$ has rank 1, when $\Lambda$ is a
maximal order, and all commutative $\Lambda$ of finite type are
maximal orders. In the graded case however, given $S$ with rank
$\qgr(S)=n$, there is no corresponding chain of graded orders
ending up with $k[X,Y]$.

We point out that if $\Cscr$ is a hereditary abelian $k$-category
with all objects of finite length, then $\Cscr$ does not
necessarily have almost split sequences (or Serre functor). For
example, this is the case if $\Cscr$ is the category of holonomic
modules over the first Weyl algebra (see \cite{P}). And it follows
from \cite{S} that it holds for the category of finite dimensional
representations over $k$ of a finite connected quiver having
oriented cycles, but which is not equal to a single oriented
cycle.

\appendix
\section{Hereditary orders and weighted projective lines}
\setcounter{section}{1}
\renewcommand{\thesection}{\Alph{section}}
\setcounter{lemma}{0} \setcounter{equation}{0}

In this appendix we will show directly that $\coh(\Oscr)$ for
$\Oscr$ a classical hereditary order over $\PP^1$ is equivalent to
$\coh\Xscr$ for a weighted projective line $\Xscr$ and furthermore
we will show that every weighted projective line appears  in this
way.  As was said before this can be deduced from Theorem
\ref{yetanothermainresult} together with
 \cite{L2}.

We follows the methods of \cite{AZ}, except that we consider
gradings by rank one abelian groups which can have torsion.

To formalize this let $\Dscr$ be an abelian category and let
$O\in\Dscr$ be an object. In addition let $(t_h)_{h\in H}$ be a
family of autoequivalences of $\Dscr$ indexed by a group $H$ and
for any pair $h_1,h_2\in H$ assume there are given natural
isomorphisms $\eta_{h_1,h_2}:t_{h_1} t_{h_2}\r t_{h_1 h_2}$
satisfying the cocycle condition
\begin{equation}
\label{cocyclecondition} \eta_{h_1h_2,h_3}(\eta_{h_1,h_2}\cdot
t_{h_3})=(t_{h_1}\cdot \eta_{h_2,h_3}) \eta_{h_1,h_2 h_3}
\end{equation}
 The data $(t_h)$,
$(\eta_{h_1,h_2})$   can be used to put a $H$-graded ring
structure on
\[
\Gamma^\ast(O)=\bigoplus_{h\in H}\Hom(t_h^{-1} O,O)
\]
as well as a $H$-graded $\Gamma^\ast(O)$-module structure on
\[
\Gamma^\ast(M)=\bigoplus_{h\in H}\Hom(t_h^{-1} O,M)
\]
>>From now on we will assume that $H$ is a finitely generated
abelian group of rank one. We fix an element $z$ in $H$.
Associated to $H$ there is a surjective map $\phi:H\r \Z$, unique
up to sign. We fix the sign by imposing  $\phi(z)>0$. If $U$ is a
$H$-graded abelian group then we say that $U$ has right bounded
grading if $U_h=0$ for $\phi(h)\gg 0$. If $R$ is a noetherian
$H$-graded ring then we define $\qgr(R)=\gr(R)/\tors(R)$ where
$\tors(R)$ consists of the right bounded modules.

The following result is an easy extension of \cite[Thm 4.5]{AZ}.
\begin{proposition} Let $\Dscr$ be a noetherian $\Ext$-finite abelian category and
  let $O\in \Dscr$. Let $(t_h)_{h\in H}$ be a system of
  autoequivalences as above and assume that  $(\Dscr,O,t_z)$ is an
  ample triple in the sense of \cite{AZ}. Then $R=\Gamma^\ast(O)$ is noetherian and the functor
  $\Gamma^\ast$ defines an equivalence between $\Dscr$ and $\qgr(R)$.
\end{proposition}
Now let $X=\PP^1$, $K=k(X)$ and let $\Oscr$ be a sheaf of
hereditary $\Oscr_X$-orders in $A=M_n(K)$.

Let $x_1,\ldots,x_t\in X$ be the set of ramification points of
$\Oscr$. Since the analysis of the cases $t=0$, $t=1$ and $t\ge 2$
is somewhat different we consider the case $t\ge 2$ first.
Afterwards we discuss the other cases.  Fix an arbitrary point $x$
in $X$ distinct from $x_1,\ldots,x_t$ and let
$(f_i)_{i=1,\ldots,t}$ be rational functions with divisor $-(x_i)+
(x)$.

Let $I_i$ be fractional $\Oscr$-ideals in $A$ defined by the
condition
\[
(I_i)_y=
\begin{cases}
(\rad \Oscr_{x_i})^{-1}&\text{if $y=x_i$}\\
\Oscr_y&\text{otherwise}
\end{cases}
\]
>>From this definition we obtain canonical isomorphisms (as
fractional ideals)
\begin{equation}
\label{canisos} I_i^{e_i}\r I_j^{e_j}:x\r xf_j/f_i
\end{equation}
We let $H$ be the abelian group of rank one generated by the
elements $h_1,\ldots, h_t$, subject to the relations
$e_ih_i=e_jh_j$ and we put $z=h_1+\cdots+h_t$.

Every $h\in H$ has a unique representation of the form
$a_1h_1+\cdots+ a_t h_t$ with $0\le a_i<e_i$ for $i>1$. We define
$I_h=I_1^{a_1}\cdots I_t^{ a_t}$. From \eqref{canisos} we obtain
canonical isomorphisms $\zeta_{h_1,h_2} :I_{h_1} I_{h_2}\r
I_{h_1+h_2}$.

Associated to the fractional ideals $I_h$ there are
autoequivalences $t_h$ on $\coh(\Oscr)$ given by $I_h\otimes-$.
The $\zeta_{h_1,h_2}$ define natural isomorphisms
$\eta_{h_1,h_2}:t_{h_1}t_{h_2}\r t_{h_1+h_2}$ satisfying the
cocycle condition \eqref{cocyclecondition}.

Now let $\bar{\Oscr}$ be a maximal order overlying $\Oscr$. As
usual $\bar{\Oscr}=\HEnd(\Escr)$ for some vector bundle $\Escr$ on
$X$. It follows from \cite[Ch IV]{RV2} that the triple
$(\coh(\Oscr),\Escr,t_z)$
is ample. Hence if we take $O=\Escr$ in the above notations and we
put $R=\Gamma^\ast(\Escr)$ then $R$ is a noetherian $H$-graded
ring and $\Gamma^\ast$ defines an equivalence between
$\coh(\Oscr)$ and $\qgr(R)$.

Our next aim will be to show that $R$ is in fact a weighted
projective line. Unfortunately the autoequivalences
$\eta_{h_1,h_2}$ clutter up our computations rather badly.
Therefore we will first give a more elegant description of $R$.

Let $D$ be the graded ring defined by
\[
D=A[u_1,u_1^{-1},\ldots,u_t,u_t^{-1}]/(f_i u_i^{e_i}=f_j
u_j^{e_j})]
\]
$D$ is clearly $H$-graded by putting $\deg u_i=h_i$.

Let $\Ascr$ be the graded order in $D$ defined by
\[
\Ascr=\bigoplus_{p_1,\ldots,p_t\in \Z} (I_1u_1)^{p_1}\cdots (I_t
u_t)^{p_t}
\]
Now it is not hard to see that
\[
R=\Hom(\Escr,\Ascr\otimes_{\Oscr}\Escr)
\]
We will determine the structure of $R$ explicitly. A local
computation shows that  $R$ is equal to
\[
\bigoplus_{p_1,\ldots,p_t\in \NN}
\Gamma(X,\Oscr_X([p_1/e_1]x_1+\cdots +[p_t/e_t] x_t))
u_1^{p_1}\cdots u_t^{p_t}
\]
where $[\alpha]$ denote the biggest integer not bigger than
$\alpha$.

We first claim that $R$ is generated by $u_1,\ldots,u_t$. By the
relations in $D$ it follows that {\small
\begin{multline*}
\Oscr_X([p_1/e_1]x_1+\cdots+[(p_i+e_i)/e_i]x_i+\cdots
+[p_j/e_j]x_j + \cdots+[p_t/e_t]x_t) u_1^{p_1}\cdots
u_i^{p_i+e_i}\cdots u_j^{p_j}\cdots u_t^{p_t}\\
=\Oscr_X([p_1/e_1]x_1+\cdots+[p_i/e_i]x_i+\cdots+[(p_j+e_j)/e_j]x_j+\cdots+[p_t/e_t]x_t)u_1^{p_1}\cdots
u_i^{p_i}\cdots u_j^{p_j+e_j}\cdots u_t^{p_t}
\end{multline*}
}%
as subsheaves of $D$. Hence  to show that every section of
$\Oscr_X([p_1/e_1]x_1+\cdots +[p_t/e_t]x_t])u_1^{p_1}\cdots
u_t^{p_t}$ is a linear combination of products of the $u_i$'s, it
suffices to do so in the case that $p_i<e_i$ for $i\ge 2$. So
below we make this assumption.

Write $p_1=q_1+ae_1$ where $q_1<e_1$.  We then have
$
\Oscr_X([p_1/e_1]x_1+\cdots +[p_t/e_t]x_t])=\Oscr_X(ax_1) $. Let
$b_1,\cdots,b_t\in \NN$ be such that $b_1+\cdots+b_t=a$. Using the
relations in $D$ we find that
\[
u_1^{q_1+b_1e_1} u_2^{p_2+b_2e_2}\cdots u_t^{p_t+b_te_t}
=
f_1^{b_2+\cdots+b_t}/(f_2^{b_2}\cdots f_t^{b_t}) u_1^{p_1}
u_2^{p_2}\cdots u_t^{p_t}
\]
The divisor of $f_1^{b_2+\cdots+b_t}/(f_2^{b_2}\cdots f_t^{b_t})$
is equal to $-a(x_1)+\sum_{i=1}^t b_i (x_i)$. Hence these rational
functions clearly generate the global sections of
$
\Oscr_X(ax_1) $, which is what we had to show.

We now claim that up to changing $f_1,(f_i)_{i\ge 3}$ by a scalar
we have the following relations in $R$~:
\begin{equation}
\label{ref:14a} u_i^{e_i}-u_2^{e_2}+\lambda_i u_1^{e_1}=0\qquad
(i\ge 3)
\end{equation}
where the $(\lambda_i)_{i\ge 3}$ are suitable scalars with
$\lambda_3=1$.

Rewriting $u_2^{e_2}$ and $u_i^{e_i}$ in terms of $u_1$ it follows
that the relation \eqref{ref:14a} is equivalent to the existence
of a linear dependence
\begin{equation}
\label{ref:15a} f_1/f_i-f_1/f_2+\lambda_i =0
\end{equation}
Now the divisors of $f_1/f_i$, $f_1/f_2$ and $1$ are respectively
given by $-(x_1)+(x_i)$, $-(x_1)+(x_2)$ and $0$. In particular
these three rational functions are all sections of $\Oscr_X(x_1)$.
Since $\Oscr_X(x_1)$ has degree one, it follows that there has to be
at least a linear dependence
\begin{equation}
\label{ref:16a} \alpha f_1/f_i+\beta f_1/f_2+\gamma=0
\end{equation}
 Furthermore, inspecting divisors, it is easily seen
that $\alpha,\beta,\gamma$ must all be non-zero.  Dividing
\eqref{ref:16a} by $-\beta$ and changing $f_i$  by a scalar yields
\eqref{ref:15a}. To make $\lambda_3$ equal to $1$ we finish by
changing $f_1$ by a suitable scalar.

At this point we know that $R$ is a quotient of the ``weighted
projective line''
\begin{equation}
\label{ref:17a} k[u_1,\ldots,u_t]/(u_i^{e_i}-u_2^{e_2}+\lambda_i
u_1^{e_1})
\end{equation}
However a straightforward computation reveals that $R$ and the
ring defined by \eqref{ref:17a} have the same Hilbert series.
Hence they are isomorphic. This concludes our analysis of the case
$t\ge 2$.

We will now discuss the other cases. First let $t=1$. We define
$R$ as above. Now $\dim R_i=1$ for $i<e_1$ and $\dim R_e=2$. Let
$v\in R_2\setminus k u_1^e$. We leave it as an exercise to the
reader to check that $R\cong k[u_1,v]$. Hence $\coh(\Oscr)$ is
again described by a weighted projective line.

The case $t=0$ is even more trivial. In that case $\Oscr$ is
Morita equivalent to $\Oscr_{\PP^1}$. So $\coh(\Oscr)$ is in fact
described by the ordinary projective line!

To finish we show that one can get all weighted projective lines
from hereditary orders. It suffices to do this in the case $t>2$.
It is convenient to choose an affine coordinate system on $\PP^1$
in such a way that $x_1=\infty$, $x_2=0$, $x_3=1$. Then up to a
scalar we have
\[
f_1(z)=z-x,\qquad f_i(z)=\frac{z-x}{z-x_i}\text{ for $i>1$}
\]
Computing the $\lambda_i$ explicitly with the above procedure we
find $\lambda_i=x_i$. This shows what we want.

\section{Examples of hereditary abelian categories}

In this appendix we give some sources of examples of hereditary
abelian categories, which are usually not $\Ext$-finite. These are
inspired by \cite{PS}.

Let $R$ be a noetherian ring of Krull dimension $n \geq 0$,
finitely generated as a module over a central subring $C$. Denote
by $\mo R$ the category of $R$-modules and as before by $\m R$ the
subcategory of finitely generated $R$-modules.  For $i \geq -1$,
let $\mathcal{C}_i$ be the subcategory of $\m R$ whose objects
have Krull dimension at most $i$, and let $\tilde{\mathcal{C}}_i$
be the subcategory of $\mo R$ whose objects are direct limits of
objects in $\mathcal{C}_i$. (We define
$\mathcal{C}_{-1}=\tilde{\mathcal{C}}_{-1}=(0)$.) Let
$\qmod_i(R)=\m R/\mathcal{C}_i$ and $\QMod_i(R)= \Mod
R/\tilde{\mathcal{C}}_i$ be the corresponding quotient categories
in the sense of \cite{G}. These are abelian categories. Similarly
we consider the case when $S$ is a \Z-graded noetherian ring
finitely generated over $k$ and finitely generated as a module
over a central subring $C$, such that $S_i=0$ for $i$ small
enough. We make the similar definitions starting with the category
$\Gr S$ of graded $S$-modules with degree zero homomorphisms, and
the subcategory $\gr S$ of finitely generated modules. The
corresponding quotient categories will be denoted by $\QGr_i(S)$
and $\qgr_i(S)$. \emph{Below when we work in the graded
  case all objects will be implicitly considered to be graded, unless
  otherwise specified.}

We can now prove the following, which gives some classes of
hereditary categories.
\begin{proposition}
\label{ref:2.5a}
   Let $R$ be a noetherian ring of Krull dimension $n \geq 0$  with
   the above assumptions and notation. Assume that $C$ satisfies
   $\Kdim C/P+\height P=n$ for every prime ideal $P$ in $C$.
Then
   the following conditions are equivalent.
   \begin{itemize}
      \item[(a)] $\QMod_i(R)$ is nonzero hereditary.
      \item[(b)] $\qmod_i(R)$ is nonzero hereditary.
      \item[(c)] Either $i=n-2$ and $\gd R_P \leq 1$ for any prime
      ideal $P$ in $C$ of height at most 1 or $i=n-1$ and $\gd R_P \leq 1$
      for any prime ideal $P$ in $C$ of height 0.
    \end{itemize}
\end{proposition}

\begin{proof}
That(a) and (b) are equivalent follows from \cite[Proposition
A3]{RV2}. To prove the other equivalences we first review some
generalities. First of all if $M$ is an $R$-module, then by
\cite[p430, Cor.\ 2]{JPGabriel} the Krull dimension of $M$ as
$R$-module is equal to the Krull dimension of $M$ as $C$-module.

Furthermore we claim that $M\in \tilde{\Cscr}_i$ if and only if
$M_P=0$ for all $P\in \Spec C$ such that $\Kdim C/P<i$ (or
equivalently, if and only if $M_P=0$ for all $P$ such that
$\height P>n-i$). To see this we may assume that $M$ is finitely
generated.  Then it follows from the theory of associated primes
that $M$ has a finite filtration (as $C$-module) with subquotients
of the form $C/Q$ with $Q\in \Spec C$. The claim is now an
immediate verification.

The subcategory $\tilde{\mathcal{C}}_i$ of $\mo R$ is a localizing
subcategory. Since $R$ is a noetherian ring which is finitely
generated as a module over its center, $\tilde{\mathcal{C}}_i$ is
closed under injective envelopes \cite[p. 431]{JPGabriel}.
Denoting by $T \colon \mo R \ra \mo R/\tilde{C}_i$ the associated
quotient functor we have that $T$ preserves injective objects and
injective envelopes by \cite[Proposition A4]{RV2}. So if $0\ra M
\ra I_0\ra I_1\ra \cdots$ is a minimal injective resolution in
$\mo R$, then $0\ra T(M)\ra T(I_0)\ra T(I_1)\ra \cdots$ is a
minimal injective resolution in $\tilde{\h_i}$. A similar
reasoning shows that for each prime ideal $P$ in $C$, we have that
$0\ra M_P\ra (I_0)_P\ra (I_1)_P\ra\cdots$ is a minimal injective
resolution in $\mo R_P$.

(c) $\Rightarrow$ (a). Assume that (c) holds, and consider for $M$
in $\mo R$ a minimal injective resolution $0\ra M \ra I_0\ra I_1
\ra \cdots$ in $\mo R$, and the induced minimal injective
resolution $0\ra M_P\ra (I_0)_P\ra (I_1)_P \ra \cdots$ for a prime
ideal $P$ in $C$.

If $i=n-2$ and $\gd R_P\leq 1$ for $\HEIGHT P \leq 1$, we get
$(I_j)_P=0$ for $j\geq 2$ and $\HEIGHT P\leq 1$, and hence $I_j$
is in $\tilde{\mathcal{C}}_{n-2}$ for $j\geq 2$, so that
$T(I_j)=0$. Then we have an injective resolution  $0 \ra T(M) \ra
T(I_0) \ra T(I_1) \ra 0$ in $\QMod_{n-2}(R)$, which shows
$\id_{\QMod_{n-2}(R)} T(M)\leq 1$. Since $T$ is essentially
surjective, it follows that $\gd {\QMod }_{n-2}(R) \leq 1$.

If $i=n-1$ and $\gd R_P \leq 1$ for $\HEIGHT P=0$, we get that
$I_j$ is in $\tilde{\mathcal{C}}_{n-1}$ for $j\geq 2$, so that $0
\ra T(M) \ra T(I_0) \ra T(I_1) \ra 0$ is an injective resolution
of $T(M)$ in $\QMod_{n-1}(R)$. Hence we have $\gd\QMod_{n-1}(R)
\leq 1$.

(a) $\Rightarrow$ (c). Assume now that $\QMod_{i}(R)$ is
hereditary, and let $0 \ra M \ra I_0 \ra I_1 \ra \cdots$ be a
minimal injective resolution of some $M$ in $\mo R$. Since then $0
\ra T(M) \ra T(I_0) \ra T(I_1) \ra \cdots$ is a minimal injective
resolution in $\QMod_i(R)$, it follows that $T(I_j)=0$ for $j>1$.
Then $I_j$ is in $\tilde{\mathcal{C}}_i$ for $j>1$, so that
 we have $(I_j)_P=0$ when $\HEIGHT P <n-i$. Thus we
 we have the exact sequence $0 \ra M_P \ra(I_0)_P \ra (I_1)_P
\ra 0$ in this case, and hence $\gd R_P \leq 1$. Since for a
noetherian ring $R$ which is finitely generated as a module over
its center we have  $\Kdim R_P \leq \gd R_P$ \cite{Brown}, and
furthermore trivially $\Kdim C_P\le \Kdim R_P$, it  follows that
$n-i-1 \leq 1$, so that $n-i \leq 2$. Hence when $\QMod_i(R) \neq
0$, we must have $i=n-2$ or $i=n-1$. In the first case we have
$\gd R_P \leq 1$ for $\HEIGHT P \leq 1$, and in the second case
$\gd R_P \leq 1$ for $\HEIGHT P=0$.
\end{proof}

In the case of graded rings $S$ with the assumptions listed in the
above, we consider graded prime ideals $P$ in the center $C$ and
graded localizations $S_P$ and graded global dimension.  Using
that also in this case $\tilde{\mathcal{C}}_i$ is closed under
injective envelopes \cite{CN}, the proof of Proposition
\ref{ref:2.5a} is easily adapted to give the following.

\begin{proposition}
\label{ref:2.6a}
   Let $S$ be a $Z$-graded ring of Krull dimension $n\geq 0$
   satisfying the standard assumptions, and with the previous
   notation. Suppose that $C$ satisfies $\Kdim C/P+\height P=n$ for
   every graded prime ideal $P$.
Then the following are equivalent.
   \begin{itemize}
      \item[(a)] $\QGr_i(S)$ is nonzero hereditary.
      \item[(b)] $\qgr_i(S)$ is nonzero hereditary.
      \item[(c)] Either $i=n-2$ and $\gd S_P\leq 1$ for any graded
        prime ideal $P$ in $C$ of height at most 1 or $i=n-1$ and $\gd
        S_P \leq 1$ for any graded prime ideal $P$ in $C$ of height 0.
   \end{itemize}
\end{proposition}
We also state the following special case.
\begin{corollary}\label{corIII.2.7}
 Let $S=C$ be a $Z$-graded commutative domain of Krull dimension 2
 satisfying the standard assumptions. Then $\qgr(S)$ is hereditary
 if and only if $S$ is an isolated singularity.
\end{corollary}

\ifx\undefined\bysame
\newcommand{\bysame}{\leavevmode\hbox to3em{\hrulefill}\,}
\fi

\end{document}